\pdfoutput=1
\RequirePackage{ifpdf}
\ifpdf 
\documentclass[pdftex]{sigma}
\else
\documentclass{sigma}
\fi

\usepackage{xspace}

\newcommand{\extcurves}{\texttt{extcurves}\xspace}

  \numberwithin{equation}{section}

\def\Im{\mathop{\rm Im}\nolimits}

\def\Jac{\mathop{\rm Jac}\nolimits}
\def\Aut{\mathop{\rm Aut}\nolimits}
\def\diag{\mathop{\rm Diag}\nolimits}

\def\id{\mathop{\rm Id}\nolimits}

\newcommand{\mi}{\mathrm{i}}
\newcommand{\me}{\mathrm{e}}
\newcommand{\dif}{\mathrm{d}}
\renewcommand{\vec}[1]{\mathbf{#1}}

\newcommand{\cx}{\bar{x}}
\newcommand{\cy}{\bar{y}}
\newcommand{\cz}{\bar{z}}

\begin{document}

\allowdisplaybreaks

\renewcommand{\thefootnote}{$\star$}

\renewcommand{\PaperNumber}{065}

\FirstPageHeading

\ShortArticleName{Bring's Curve: its Period Matrix and the Vector of Riemann Constants}

\ArticleName{Bring's Curve: its Period Matrix\\ and the Vector of Riemann Constants\footnote{This
paper is a contribution to the Special Issue ``Geometrical Methods in Mathematical Physics''. The full collection is available at \href{http://www.emis.de/journals/SIGMA/GMMP2012.html}{http://www.emis.de/journals/SIGMA/GMMP2012.html}}}

\Author{Harry W.~BRADEN and Timothy P.~NORTHOVER}
\AuthorNameForHeading{H.W.~Braden and T.P.~Northover}
\Address{School of Mathematics, Edinburgh University, Edinburgh, Scotland, UK}
\Email{\href{mailto:hwb@ed.ac.uk}{hwb@ed.ac.uk}, \href{mailto:T.P.Northover@gmail.com}{T.P.Northover@gmail.com}}

\ArticleDates{Received June 10, 2012, in f\/inal form September 27, 2012; Published online October 02, 2012}

\Abstract{Bring's curve is the genus 4 Riemann surface with automorphism group
of maximal size, $S_5$. Riera and Rodr{\'{\i}}guez have provided the most detailed study
of the curve thus far via a hyperbolic model. We will recover and extend their results via
an algebraic model based on a sextic curve given by both Hulek and Craig and implicit in work of Ramanujan.
In particular we recover their period matrix; further,
the vector of Riemann constants will be identif\/ied.}


\Classification{14H45; 14H55; 14Q05}

\renewcommand{\thefootnote}{\arabic{footnote}}
\setcounter{footnote}{0}

\section{Introduction}

Bring's curve is the genus 4 Riemann surface with the automorphism group
of maximal size, $S_5$ \cite{breuer, edge78, edge81, weber05}. It may be expressed as the
complete intersection in $\mathbb{P}^4$ given by
\begin{gather*}
  x_1 + x_2 + x_3 + x_4 + x_5  = 0, \\
  x_1^2 + x_2^2 + x_3^2 + x_4^2 + x_5^2  = 0, \\
  x_1^3 + x_2^3 + x_3^3 + x_4^3 + x_5^3  = 0.
\end{gather*}
Here the permutations of the coordinates
$x_i$ make manifest the $S_5$ symmetry. The curve naturally arises in the study of the
general quintic $\prod\limits_{i=1}\sp5(x-x_i)$ when this is reduced to Bring--Jerrard form $x^5+bx+c$.
Just as with Klein's curve, Bring's curve
may be studied by either plane algebraic or hyperbolic models. Perhaps the most detailed study thus far
is that of Riera and Rodr{\'{\i}}guez~\cite{riera} via a hyperbolic model. Using a
representation much like that of Klein's curve they produced the very
simple period matrix
\begin{equation}\label{RRstau}
  \tau = \tau_0
  \begin{pmatrix}
    4 & 1 & -1 & 1 \\
    1 & 4 & 1 & -1 \\
    -1 & 1 & 4 & 1 \\
    1 & -1 & 1 & 4
  \end{pmatrix},
\end{equation}
for a determined $\tau_0\in\mathbb{C}$. This period matrix already
exhibits much of the symmetry implicit in the automorphism group and
we won't attempt to improve on this result. We shall however reproduce this result
and the homology basis of Riera and Rodr{\'{\i}}guez by studying a plane model of the curve
and then compute the vector of Riemann constants.
All these we believe are new. To do this we shall use and extend the techniques of \cite{bn10}.
These techniques have been developed to implement the modern approach to integrable systems
based upon a spectral curve.

It remains to introduce the plane model of Bring's curve we shall employ. In \cite{dye}, Dye
explicitly gives a sextic plane curve and proves its equivalence to Bring's.
The remarkable fact about this representation is that, of the full $S_5$ symmetry group, $A_5$ is
generated by projectivities in $\mathbb{P}^2$. Dye's sextic is not the one we use.
In \cite{craig}\footnote{See~\cite{craig04} for errata; see
 \url{http://members.optusnet.com.au/~towenaar/} for a corrected version.}, Craig studies the rational points of a second genus-4
sextic which possesses at least $A_5$ as a symmetry group. This work generalizes a similar
result for Klein's curve, where the curve is parameterized by modular functions. Craig
observes that work of Ramanujan means the coordinates of the curve may also be expressed in terms
of modular functions. In fact
Craig's model is very closely related to Dye's and we will
show that it too is equivalent to Bring's curve by giving an explicit
transformation of $\mathbb{P}^2$ mapping between the representations of Dye and Craig.
The sextic studied by Craig had in fact been introduced by Hulek \cite[p.~82]{hulek85}
who also makes connection with the modular properties, and we will refer to this curve as the Hulek--Craig
(HC) curve throughout. This representation will be more useful for our purposes than Dye's
since it has a more obvious real structure and simpler branching
properties.

An outline of the paper is as follows: in Section~\ref{section2} we shall discuss some plane sextics describing Bring's curve.
The Hulek--Craig curve will be described in detail in Section~\ref{section3} while in Section~\ref{section4} we shall recall the Riera--Rodr\'iguez
hyperbolic model of Bring's curve. Here a detailed analysis will enable us to identify the two descriptions
and in particular the homology basis of Riera and Rodr\'iguez, the period matrix then following. Our identif\/ication will make use of the real structure and f\/ixed oval of the models, described in increasing detail in Sections~\ref{section2} and~\ref{section3}. Finally in
Section~\ref{section5} we determine the vector of Riemann constants for the curve.

\section{Two sextics} \label{section2}

\subsection{Dye's sextic}\label{section2.1}

Let $j = \frac{1+\sqrt{5}}{2}$, a root of $j^2=j+1$. Dye \cite{dye} introduces
 the plane sextic curves given by
\begin{gather*}
  \mathcal{D}_\lambda(x,y,z) :=
  (x+jy)^6 + (x-jy)^6 + (y+jz)^6  \\
\hphantom{\mathcal{D}_\lambda(x,y,z) :=}{}    + (y-jz)^6 + (z+jx)^6 + (z-jx)^6 + \lambda\big(x^2+y^2+z^2\big)^3 = 0.
\end{gather*}
For generic $\lambda\in\mathbb{C}$ the curve has genus 10, but if
$\lambda$ is chosen to be $-\frac{78+104j}{5}$ then the genus drops to
4 and the resulting curve is shown to be equivalent to Bring's. We
correspondingly def\/ine
\[
  \mathcal{D}(x,y,z) := \mathcal{D}_{-\frac{78+104j}{5}}(x,y,z).
\]
The curve $\mathcal{D}(x,y,z)$ has the obvious order three cyclic symmetry
\[ b' : \ (x,y,z) \mapsto (y,z,x), \]
as well as the less obvious order two symmetry
\[ a' : \
  \begin{pmatrix}
    x\\y\\z
  \end{pmatrix}
  \mapsto
  \begin{pmatrix}
    -j  & 1    & j^2 \\
    1   & -j^2 & j \\
    j^2 & j    & 1
  \end{pmatrix}
  \begin{pmatrix}
    x\\y\\z
  \end{pmatrix},
\]
both presented by Dye in his paper. It is easy to check that these are the
classical generators for $A_5$:
\begin{equation*}
  a'b' =
  \begin{pmatrix}
    j^2 & -j  & 1 \\
    j   & 1   & -j^2 \\
    1   & j^2 & j
  \end{pmatrix}
\end{equation*}
has order f\/ive and (taking into account the projective nature of the
space)
\[ a'^2 = b'^3 = (a'b')^5 = 1. \]

\subsection[The Hulek--Craig sextic]{The Hulek--Craig sextic}\label{section2.2}

Hulek and Craig both introduce the sextic
\begin{equation}
\label{craig}
  \mathcal{C}(\cx,\cy,\cz) :=
    \cx\big(\cy^5+\cz^5\big) + (\cx\cy\cz)^2 - \cx^4\cy\cz -2(\cy\cz)^3  = 0.
\end{equation}
This curve is also of genus 4
and admits~$A_5$ as a symmetry group\footnote{Here a bar over
  variables is used to distinguish them from the Dye curve, rather
  than to denote complex conjugation. Craig notes the results of Ramanujan \cite[Chapter~19, Entry 10(iv), 10(vii)]{berndtiii}
entail the parameterization
\begin{gather*}
(\cx,\cy,\cz)
 =
\left(\sum_{n=-\infty}\sp{\infty} q^{(5n)^2},\sum_{n=-\infty}\sp{\infty} q^{(5n+1)^2},
\sum_{n=-\infty}\sp{\infty} q^{(5n+2)^2}\right).
\end{gather*}}.
In this case an order f\/ive symmetry is obvious and we may take
\begin{equation*}
  \overline{ab} : \ (\cx,\cy,\cz) \mapsto \big(\zeta^2 \cx, \zeta^4 \cy, \cz\big),
\end{equation*}
where $\zeta=\me^{2\pi\mi/5}$.
There is also a corresponding order two symmetry
\[ \overline{a} : \
  \begin{pmatrix}
    \cx\\\cy\\\cz
  \end{pmatrix}
  \mapsto
  \begin{pmatrix}
    1 & 2                 & 2 \\
    1 & \zeta+\zeta^{-1}   & \zeta^2+\zeta^{-2} \\
    1 & \zeta^2+\zeta^{-2} & \zeta+\zeta^{-1}
  \end{pmatrix}
  \begin{pmatrix}
    \cx\\\cy\\\cz
  \end{pmatrix}.
\]
Together these generate $A_5$ again since $\overline{a}\overline{ab} =:
\overline{b}$ has order three (hence the slightly unusual choice of
notation for $\overline{ab}$ above).

In fact the sextic (\ref{craig}) is also a model for Bring's curve using Dye's result and the following theorem.
\begin{theorem} With
$
  A =
  \begin{pmatrix}
    j & 1              & 1 \\
    0 & -\mi\sqrt{2+j} & \mi\sqrt{2+j} \\
    1 & -j             & -j
  \end{pmatrix}
$
then
$
  \mathcal{D}(A\vec{x}) = -960(9 + 4\sqrt{5}) \mathcal{C}(\vec{x})$
{and hence}
  $\mathcal{D}(\vec{x}) = 0 \iff \mathcal{C}(A^{-1} \vec{x}) = 0
$.
\end{theorem}
This may be directly verif\/ied. The choice of the matrix $A$ follows upon consideration of the
conjugacy classes of automorphisms of both models (see \cite{northover} for further details).

We observe that the antiholomorphic involution $[\cx,\cy,\cz]\mapsto [\cx\sp\ast,\cy\sp\ast,\cz\sp\ast]$
(where $\sp\ast$ is complex conjugation) is a symmetry of $\mathcal{C}$ though it is orientation-reversing and
so not a conformal automorphism.  This involution endows $\mathcal{C}$ with a real structure. The f\/ixed
 point set of such a real structure is either empty or the disjoint union of simple closed curves, known as ovals
 following Hilbert's terminology \cite{begg}. A classical result of Harnack for a Riemann surface of genus $g$ with real structure says there are at most $g+1$ ovals. We shall show the the HK curve has one oval with this real structure.

\section[Details of the Hulek-Craig curve]{Details of the Hulek--Craig curve}\label{section3}

The representation (\ref{craig}) will turn out to be the most convenient for later
work so it is worth spending some time on its detail, particularly its
desingularisation.

\subsection{Special points of the Hulek--Craig representation and
desingularisation}\label{section3.1}

The points at inf\/inity for the HC curve \eqref{craig} are
given by (the real points) $[0,1,0]$ and $[1,0,0]$, but the latter is singular.  In
fact the singularities of the HC curve are $[1,0,0]$ and $[\zeta^k,
\zeta^{2k}, 1]$ for $k\in\{0,\dots,4\}$ so we must work out expansions
nearby in order to form a properly compact curve.

First the inf\/inite singularity: consider the structure near $[1,0,0]$,
say at points $[1,y,z]$ for small $y$, $z$. The curve reduces to
\[
  y^5 + z^5 + y^2 z^2 - yz -2y^3 z^3 = 0,
\]
so in the usual Puiseux construction we suppose $z = A y^{\alpha} +
\cdots$. Equating lowest order terms we get one of
\begin{itemize}\itemsep=0pt
\item $A^5 y^{5\alpha} - Ay^{\alpha+1} = 0$ which implies $z =
  y^{1/4}+\cdots$, that is $z\approx y^{1/4}$ near this point,
\item $y^5 - Ay^{\alpha+1} = 0$ which implies $z = y^4 + \cdots$ or $z\approx y^{4}$.
\end{itemize}
The second of these gives a single $z$ for each $y$ near 0, the f\/irst
gives four dif\/ferent values for $z$. Together these make up the
expected f\/ive sheets and so expansions after this point are
unique. Thus in the vicinity of $[1,0,0]$ solutions $[1,y,z]$
of the f\/irst equation behave as $[1, t^4, t]$ where $t$ is a local parameter
for the curve. Similarly the second equation has solutions behaving as $[1,t,t^4]$ in terms of
a local parameter. Thus the point $[1,0,0]$ desingularises into precisely
two points on the nonsingular curve:
\begin{gather}
  \label{eq:inftyExp}
    [1,0,0]_1  \sim \big[1,t^4,t\big], \qquad
    [1,0,0]_2  \sim \big[1,t,t^4\big],
\end{gather}
where $\sim$ here indicates behaviour of a local coordinate in the vicinity of a specif\/ied point.

For the remaining singular points we only need to investigate
explicitly one and then note that the automorphism $[x,y,z] \mapsto
[\zeta x, \zeta^2 y, z]$ will tell us how the other singularities
behave. So we look at $[1,1,1]$. At f\/irst sight, two of the sheets come
together here. Consider $[1+\epsilon,y, 1]$ near to $[1,1,1]$. To
f\/irst order
\begin{equation*}
  y^5- 2y^3 + y^2 - y + 1 = 0.
\end{equation*}
This quintic has four distinct roots: two are complex, corresponding to
nonsingular points and will play no role in future developments. There
is a real negative root $\alpha\sim -1.7549$ which also corresponds
to a nonsingular point and will occur later. Finally, 1 is a root,
which gives us the expected singularity at $[1,1,1]$.

Expanding about this singular point, at the next order we discover
\begin{equation*}
  y = 1 + \epsilon\frac{-1+\sqrt{5}}{2} + \cdots, \qquad
  y = 1 + \epsilon\frac{-1-\sqrt{5}}{2} + \cdots.
\end{equation*}
These are clearly distinct solutions and together with the nonsingular
expansions exhaust the f\/ive possible nonsingular preimages near $x=1$,
so $[1,1,1]$ once again desingularises to two distinct points.

\subsection[Branched covers of $\mathbb{P}^1$]{Branched covers of $\boldsymbol{\mathbb{P}^1}$}\label{section3.2}

We now consider the curve (\ref{craig}) as a branched cover of $\mathbb{P}^1$
with $x$ as the coordinate. The af\/f\/ine part of the HC curve is obtained by setting $z=1$
in \eqref{craig} yielding
\begin{equation}
  \label{craigp1}
  x y^5 + x + x^2 y^2 - x^4 y - 2 y^3 = 0.
\end{equation}
There are 5 sheets above the generic $x$, with branch points at $0,
\infty$ and
\[
  \left\{ \frac{\zeta^k}{4}
    \big(1674\pm 870\mi\sqrt{15}\big)^{1/5}
    : \;  k \in \{0,\dots,4\}
  \right\}.
\]
There is also a double solution at $x=\zeta^k$ but these are
precisely the singular points similar to $[1,1,1]$ we investigated
before and after resolution the cover is regular there.

At $x=0$ we have two preimages, one corresponding to $[0,0,1]$ with an
expansion
\begin{equation}
  \label{eq:x0exp1}
  y = \frac{1}{2^{1/3}} x^{1/3} + \cdots,
\end{equation}
where three sheets come together, and the other corresponding to $[0,1,0]$ with expansion
\begin{equation}
  \label{eq:x0exp2}
  y = \sqrt{2} x^{-1/2} + \cdots
\end{equation}
where two sheets come together.
Similarly at $x=\infty$ we have two preimages after
desingularisation: $[1,0,0]_1$ where from $[1,t^4,t]\equiv [t\sp{-1},t\sp{3},1]$ and
$y\approx x\sp{-3}$; and $[1,0,0]_2$ where from $[1,t,t^4]\equiv [t\sp{-4},t\sp{-3},1]$ and
$y\approx x\sp{3/4}$.
The other 10 branch
points correspond to the solutions of $256 x^{10}-837 x^5+3456=0$ and have two sheets
coming together at each\footnote{We remark that the Maple command {\tt{monodromy}$(x y^5 + x + x^2 y^2 - x^4 y - 2 y^3,x,y,$\tt{showpaths})} will produce the monodromy data for the HK curve, together with the paths and sheet numbering
necessary to make sense of this data. The branch point $0$ has monodromy $[1,2][3,4,5]$ while that of $\infty$ is $[1,4,5,2]$; these are the cycle structures described above. The remaining ten branch points arranged with increasing argument have monodromies $[1,4]$, $[2,4]$, $[2,5]$, $[1,5]$, $[1,3]$, $[2,3]$,$[2,4]$, $[1,4]$, $[1,5]$ respectively,
here indicating the two sheets that come together.}.

\subsection[Real paths on the Hulek-Craig curve]{Real paths on the Hulek--Craig curve}\label{section3.3}

The real structure of the HK curve leads to real ovals. Here we shall show there is in fact one.
We begin by looking at portions of this oval, which we shall refer to as a `real path' and then
indicate how they join together\footnote{The Maple command {\tt{plot\_real\_curve}}$(x y^5 + x + x^2 y^2 - x^4 y - 2 y^3,x,y)$ will in fact plot this directly.}.
The real paths in this cover will be of particular interest later on
so we will take some time to explore their nature now. We begin with the af\/f\/ine curve, further noting what happens at the real inf\/inite points $[0,1,0]$, $[1,0,0]_{1,2}$ which compactify the real curve.

First, if the
number of real roots of \eqref{craigp1} considered as a polynomial in
$y$ changes then its discriminant
\begin{gather}
 \Delta(x) = -x^3\big(256 x^{20}-1349x^{15}+5386x^{10}-7749x^5+3456\big)\nonumber\\
 \hphantom{\Delta(x)}{}  =
 -x^3\big(256 u^2-837 u+3456\big)(u-1)^2, \qquad u=x^5,  \label{eq:disc}
\end{gather}
must vanish there. The only real roots of the equation $\Delta(x)=$ are $x=0,1$,
so we are reduced to considering the intervals $(-\infty,0), (0,1),
(1,\infty)$.
\begin{itemize}\itemsep=0pt
\item If $x < 0$ then there is just one real root.
\item If $0 < x < 1$ then there are three real roots.
\item If $x > 1$ then there are also three real roots.
\end{itemize}

Referring to the expansions \eqref{eq:x0exp1} and \eqref{eq:x0exp2} we
see that a real path starting with $x < 0$ moving towards $x=0$ must
be approaching $[0,0,1]$ along the expansion $y = 2^{-1/3} x^{1/3}
+\cdots$ (i.e.\ $y\rightarrow 0$ too). Continuity demands that when
extended past $x=0$ it  too should have $y$ small and positive for small
$x>0$. We will call this path $\gamma_0$.

We now turn our attention to another real path approaching $x=0$, this
time for $x>0$. It must lie on the expansion $y=\sqrt{2} x^{-1/2} +
\cdots$ and hence $y$ is either large and positive or large and
negative; we will call these paths $\gamma_{+}$ and $\gamma_{-}$. In
fact the expansion is telling us that $\gamma_{+}$ and $\gamma_{-}$
meet at $[0,1,0]$ and we could form a single continuous path, but we
will maintain the distinction for now.

In summary we have three real paths coming out of $x=0$ along the
positive axis, satisfying (for small $x > 0$),
\[
  y(\gamma_{-}) \ll 0 < y(\gamma_0) \ll y(\gamma_+).
\]

Now we are ready to consider what happens at $x=1$. On the
desingularised curve there are three real points here (the two from
desingularising $y=1$ and the remaining real root $\alpha$).
Each of the curves coming out of $x=0$ must pass through one of them.
Further, the order of the $y$ values among the paths must be the same
approaching $x=1$ as it was leaving $x=0$ since, otherwise, they would
have crossed in between and this would have shown itself in~\eqref{eq:disc}.

The three expansions near $x=1$ in order of increasing $y$ for $x <
1$ are
\begin{gather*}
  y  \approx \alpha,\qquad
  y  \approx 1 + (x-1)\frac{-1+\sqrt{5}}{2}, \qquad
  y  \approx 1 + (x-1)\frac{-1-\sqrt{5}}{2}.
\end{gather*}
Thus the path that started $y\ll 0$ must pass through the f\/irst point,
$y\approx 0$ must pass through the second and $y\gg 0$ the third.
Signif\/icantly this means the latter two paths actually cross at $x=1$
and for $x=1+\epsilon$ we have
\[
y(\gamma_{-}) < y(\gamma_{+}) < 1 < y(\gamma_0).
\]

Finally we consider the remaining points $[1,0,0]_{1,2}$ at $\infty$. Recall the expansions~\eqref{eq:inftyExp}. If $x \ll 0$ then naturally there is only one
real path, which arrives at $[1,0,0]_1$ with small $y$. If $x \gg 0$
the situation is very similar to $x=0$: two expansions with $|y| \gg
0$ arriving at $[1,0,0]_2$ and one lying between these with $y \approx 0$. As
before, the paths cannot have crossed between $x=1$ and $x=\infty$ and
so we are forced to conclude that $\gamma_-$ has the expansion $y\approx
-x^{3/4}$, $\gamma_+$~has the expansion $y\approx x^{-3}$ and $\gamma_0$
has the expansion $y\approx x^{3/4}$ near~$\infty$.

\begin{figure}[t]
  \centering
  \includegraphics[scale=0.93]{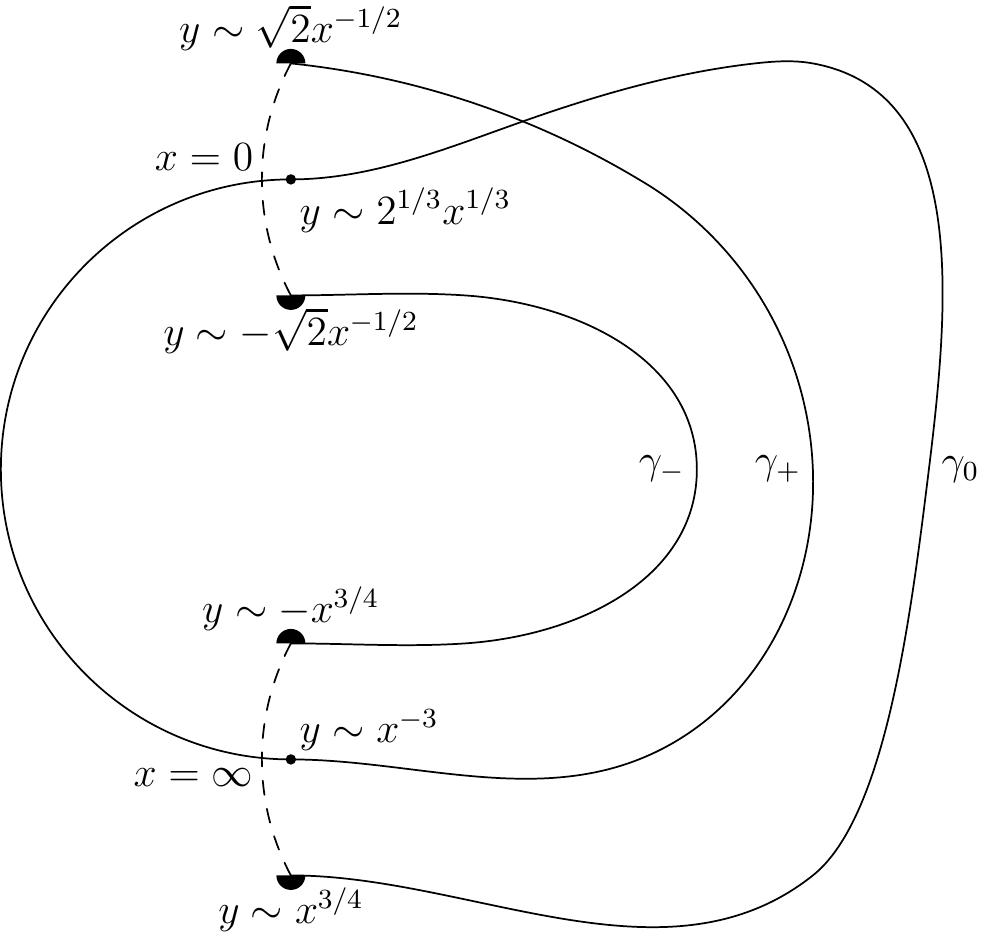}
  \caption{Real paths on the Hulek--Craig curve as a branched cover of $\mathbb{P}^1$.}
  \label{fig:craigreal}
\end{figure}

Putting these facts together we can plot Fig.~\ref{fig:craigreal}
(the joined semicircular dots represent the same point on the curve,
separated to show the distinct $y$ values of paths entering them).
We discover that all the paths ($\gamma_-$, $\gamma_0$, $\gamma_+$ and the
$x<0$ path) actually form part of one large closed loop showing that the real structure of
the HK curve has one oval. (Another proof of this will be given in the next section.)

\section{The Riera and Rodr{\'{\i}}guez hyperbolic model}\label{section4}

\subsection[Introduction to $H$]{Introduction to $\boldsymbol{H}$}\label{section4.1}

Riera and Rodr{\'{\i}}guez, in \cite{riera}, give Bring's curve as a
quotient, $H$, of the hyperbolic disc. They then proceed to calculate
a period matrix taking account of the symmetries of the curve.

\begin{figure}[t]
  \centering
  \includegraphics[scale=0.93]{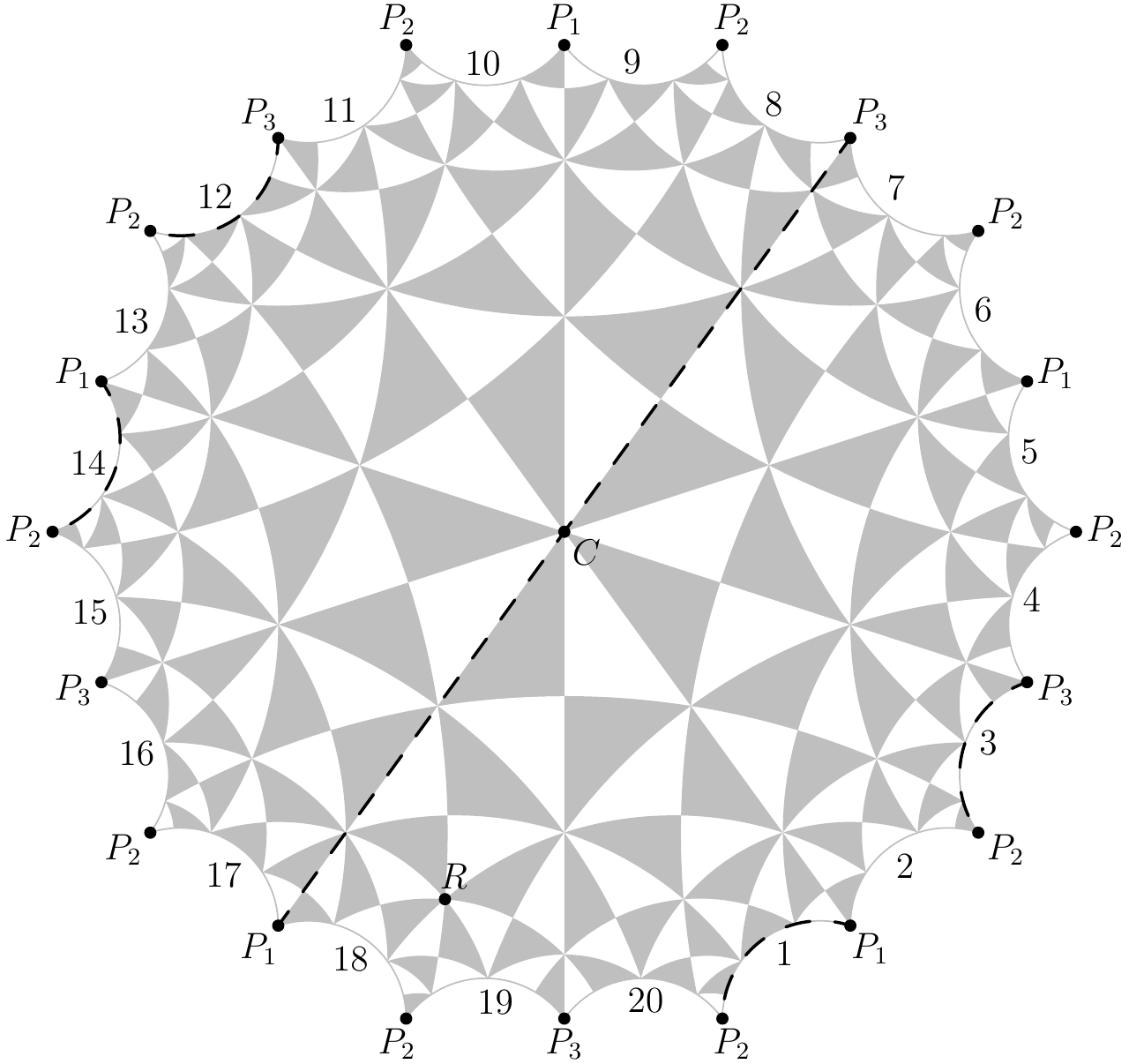}

  \begin{tabular}{rl|rl|rl|rl|rl}
    \multicolumn{10}{c}{Edge identif\/ications}\\
    \hline
    1 &$\leftrightarrow$ 14 &
    5 &$\leftrightarrow$ 18 &
    9 &$\leftrightarrow$ 2 &
    13 &$\leftrightarrow$ 6 &
    17 &$\leftrightarrow$ 10 \\
    3 &$\leftrightarrow$ 12 &
    7 &$\leftrightarrow$ 16 &
    11 &$\leftrightarrow$ 20 &
    15 &$\leftrightarrow$ 4 &
    19 &$\leftrightarrow$ 8
  \end{tabular}
  \caption{Riera and Rodr{\'{\i}}guez hyperbolic model, $H$, of Bring's curve.}
  \label{fig:bringHyperbolic}
\end{figure}

The essential features of the model can be seen in Fig.~\ref{fig:bringHyperbolic}. The surface is seen to be a 20-gon with
edges identif\/ied as shown in the table below the f\/igure. (We refer, for example, to the
identif\/ied edges~$2$ and~$9$ as~$2/9$.)  This leads to
the polygon's vertices falling into three equivalence classes, also
annotated in the f\/igure.  Naturally, this surface has genus~4.

\begin{figure}[t]
  \centering
  \includegraphics[scale=0.93]{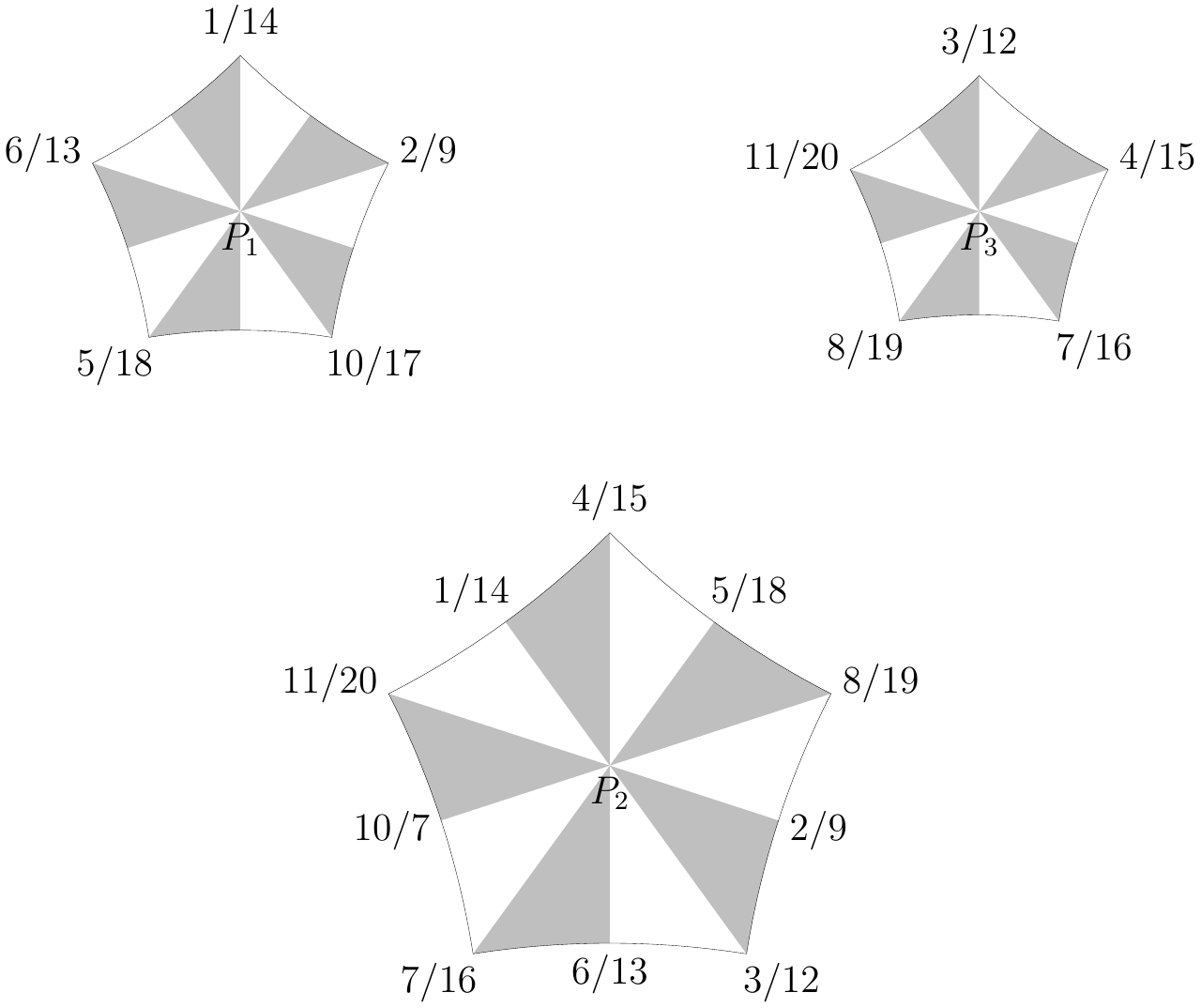}
  \caption{Conformal structure of $P_1$, $P_2$ and $P_3$.}
  \label{fig:ptreconst}
\end{figure}

For future calculations it will also be very useful to know the
conformal structure (or equiva\-lently, the local holomorphic coordinate)
about the points~$P_1$,~$P_2$ and~$P_3$. This can be
reconstructed quite easily from Fig.~\ref{fig:bringHyperbolic}. For
example, start near $P_1$ in the bottom right quadrant on edge~$2/9$. Make a small arc around $P_1$ proceeding anticlockwise and you
will next reach edge~$1/14$. Repeating at edge~14 tells us that we
next meet~$6/13$. If this procedure is continued we obtain Fig.~\ref{fig:ptreconst}.

The polygon can be tiled by 120 double triangles (one can take a sector of
the central pentagon as a fundamental domain). Now consider the
automorphism group. Let $d$ be a rotation of $\frac{\pi}{2}$ about a
vertex of the central pentagon and $c$ be a rotation of $\pi$ about
the midpoint of an adjacent pentagon edge. Then clearly $c^2 = d^4 =
1$. But it is also easy to see that $cd$ is a rotation of
$\frac{2\pi}{5}$ about the centre and hence $(cd)^5 = 1$. The rotations~$c$ and~$d$
are thus the classical generators of~$S_5$ and this describes the
entire automorphism group of Bring's curve.

Riera and Rodr{\'{\i}}guez give the homology basis for this model by
prescribing which edges of the polygon to traverse. We are going to
construct an equivalent basis for the HK curve by understanding an
isomorphism
\begin{equation}
\label{hypiso}
f : \ H \rightarrow
  \left\{(x,y,z)\in\mathbb{C}^3 : \mathcal{C}(x,y,z) = 0 \right\}
\end{equation}
well enough to determine the precise values to which each edge of the polygon
in Fig.~\ref{fig:bringHyperbolic} maps. Once this is achieved,
converting the homology basis will be a purely mechanical af\/fair as illustrated
in~\cite{bn10} for Klein's curve. Along the way we will gain some understanding of how~$f$ acts on the automorphism group by push-forwards.

\subsection{Riera and Rodr{\'{\i}}guez basis}\label{section4.2}

We start by recapitulating the hyperbolic basis of interest.
Riera and Rodr{\'{\i}}guez begin with a simple non-canonical basis. They
f\/irst def\/ine
\begin{gather*}
  \alpha_1  = 1 + 2, \qquad
  \alpha_2  = 3 + 4
\end{gather*}
(in edge traversal notation, see~\cite{rauch, riera}) and then act on these cycles by rotations of
$\frac{2\pi k}{5}$ to obtain their initial basis. So essentially
\[
  \alpha_i = (2i-1) + (2i).
\]

Next they specify (by f\/iat) a matrix which transforms these~$\alpha_i$
into a canonical basis and proceed to derive further basis change to
make use of the symmetries. The end result is the following
basis-change matrix (implicit in~\cite{riera})
\begin{equation}
  \label{basistrans}
  \begin{pmatrix}
  1 & 0  & 0 & 0  & -2 & 0 & 1  & 0  \\
  1 & -1 & 0 & -1 & -1 & 1 & 1  & 1  \\
  1 & -1 & 0 & 0  & -1 & 2 & 1  & -1 \\
  0 & -1 & 0 & 0  & 1  & 2 & 0  & 0  \\
  1 & -1 & 1 & 0  & -1 & 1 & -1 & -1 \\
  1 & -1 & 1 & -1 & 0  & 0 & -1 & 1  \\
  1 & -1 & 1 & -1 & 0  & 1 & 0  & 1  \\
  0 & -1 & 0 & -1 & 1  & 1 & 1  & 2
  \end{pmatrix},
\end{equation}
which sends the initial $\alpha_1,\dots,\alpha_{8}$ homology basis to
another $\{\mathfrak{a}_i,\mathfrak{b}_i\}_{i=1}\sp{4}$, that is not only canonical but behaves well with respect to
the symmetry group of the curve. Now the symmetries relate the periods $\mathcal{A}_{ij}=\int_{\mathfrak{a}_i}v_j$ and
$\mathcal{B}_{ij}=\int_{\mathfrak{b}_i}v_j$ (for any basis of holomorphic dif\/ferentials $v_i$).
As a~consequence, the period matrix $\tau=\mathcal{B}\mathcal{A}\sp{-1}$
can be
written as \eqref{RRstau},
where $\tau_0 \approx -0.5+0.185576\mi$ is def\/ined in terms of
Klein's $j$-invariants\footnote{Riera and Rodr{\'{\i}}guez swap these two
equations. However, we believe this to be a typographical error.} by
\begin{gather}\label{deftau0}
  j(\tau_0) = -\frac{29^3\times 5}{2^5}, \qquad
  j(5\tau_0) = -\frac{25}{2}.
\end{gather}

\subsection[Understanding the isomorphism $f$]{Understanding the isomorphism~$\boldsymbol{f}$}\label{section4.3}

We now turn our attention to the isomorphism, $f$, mentioned in~\eqref{hypiso}. Clearly  there won't be a single isomorphism
since if $a$ is an automorphism of $H$ and $\sigma$ of an automorphism of
the HC
curve then $\sigma\circ f \circ a$ will also be an
isomorphism from the hyperbolic model $H$ to the HC
representation. We will exploit this fact.

There are two key ingredients to our identif\/ication. First is the rotation of the entire
hyperbolic polygon about its centre by $2\pi/5$ (the automorphism $cd$
above). This automorphism allows us to express all twenty of the
polygon's edges in terms of just four, a great simplif\/ication of our
problem. If we knew the values of $f$ on four edges, and the matrix
representing $f_*(cd)$, the induced action of $cd$ on the HC curve, then
\begin{align*}
  f(\text{edge } k+ 4) &= f\left((cd)(\text{edge } k)\right)
    =  f_*(cd) f(\text{edge } k),
\end{align*}
which allows us to compute the values of $f$ on the remaining 16
edges.

Second is a geodesic ref\/lection of the hyperbolic disc which will correspond to our real structure;
the geodesic is denoted by the
dashed lines in Fig.~\ref{fig:bringHyperbolic}. The line starts at $P_3$,
goes through~$C$ to~$P_1$, along edge~1 to~$P_2$ and along edge 3 back
to $P_3$. If we knew how this acted on the HC model, we would know
its f\/ixed points correspond in some manner to edges~1/14 and~3/12, and
the marked diameter.
Identifying points~$P_1$, $P_2$ and~$P_3$ on the HC representation
would then complete the picture by dividing this f\/ixed line up into
just the intervals needed to draw homology paths around known branch
points.

Starting with the central rotation~$cd$ on the hyperbolic
model and some isomorphism~$f$ to the HC representation, since all
order 5 elements of $S_5$ are conjugate there is an HC-automorphism
${\sigma}\in S_5$ such that
\[
  {\sigma} f_*(cd) {\sigma}^{-1} = Z^k,
\]
where $k\in\{0,\dots,4\}$ and
\[ Z : \ [x,y,z] \mapsto \big[\zeta x, \zeta^2 y, z\big]. \]
We are being f\/lexible about which power of $Z$ occurs here because
later choices (specif\/ically rotations about $R$ in Fig.~\ref{fig:bringHyperbolic}) will modify any decision made at
this stage. But then
\begin{gather*}
  (\sigma\circ f)_*(cd)  = \sigma_*(f_*(cd))
    = \sigma f_*(cd) \sigma^{-1}
    = Z^k.
\end{gather*}
So the isomorphism $\sigma\circ f$ from the hyperbolic model to
the HC model sends $cd$ to $Z^k$.

Now consider a rotation about $R$ in Fig.~\ref{fig:bringHyperbolic} which
cyclically permutes the f\/ixed points of $cd$. The f\/ixed points
on the hyperbolic side are $C$, $P_1$, $P_2$, $P_3$ and on the HC side
$[0,1,0]$, $[0,0,1]$, $[1,0,0]_1$, $[1,0,0]_2$.  Let integers $i$ and $n$ be
def\/ined by the equations
\begin{gather*}
  P_i  = (\sigma\circ f)^{-1}([0,0,1]),  \qquad
  R^n(P_i)  = C.
\end{gather*}
Then
\begin{gather*}
  (\sigma\circ f\circ R^{-n})(C)
     = (\sigma\circ f\circ  R^{-n})(R^n(P_i))
    = (\sigma\circ f)(P_i)
    = [0,0,1],
\end{gather*}
and further
\begin{gather*}
  \left(\sigma\circ f\circ R^{-n}\right)_*(cd)
     = (\sigma\circ f)_*\left(R^{-n}_*(cd)\right)
    = (\sigma\circ f)_*\left((cd)^j\right)
    = Z^{jk}
    = Z^m,
\end{gather*}
for some integers $j$ and more importantly~$m$. Since we haven't f\/ixed
the power of $Z$ up to now this means that $\sigma\circ f\circ R^{-n}$
serves our purposes just as well as~$\sigma\circ f$ did.

Although we have used most of the available freedoms to constrain the relation bet\-ween~$f$,~$Z$ and~$C$,  we actually still have the ability to apply a central
rotation, if it would help since that would alter neither of the
properties above.

Next consider complex conjugation on the HC model. This is a symmetry that
reverses orientation (and so not part of the~$S_5$ symmetry group). It
f\/ixes an entire line (the real axis) including the f\/ixed points of~$Z$. In the hyperbolic picture this means it must be a ref\/lection
about some diameter. We use our f\/inal remaining freedom to demand that
it is ref\/lection about the dashed diameter in Fig.~\ref{fig:bringHyperbolic}, i.e.\ that the real axis in the HC model
corresponds to these dashed edges (and diameter).

We now have two tasks remaining:
\begin{itemize}\itemsep=0pt
\item Find out what $P_1$, $P_2$ and $P_3$ become on the HC model so
  we can describe edge 1/14 as the real path from $P_2$ to $P_1$ and edge
  3/12 as the real path from $P_2$ to $P_3$.
\item Find out what power of $Z$ the central rotation of $2\pi/5$
  becomes so we can describe (for example) edge 4/15 as $Z^k$ applied
  to the real path from $P_2$ to $P_3$.
\end{itemize}

The second task is actually easier to accomplish at this
stage. Consider the structure near $[0,0,1]$ (which we demanded was
the centre of the polygon, $C$, hyperbolically); there are three
sheets coming together at this branch so unwrapping it will
ef\/fectively divide angles by~3. Ma\-the\-matically this means that any set
of manifold coordinates $\phi : \mathcal{C} \rightarrow \mathbb{C}$
centred on $[0,0,1]$ will satisfy
\[
  \phi([x,y,1])^3 = \alpha x + O\big(x^2\big).
\]
In these coordinates, since $[0,0,1]$ is a f\/ixed point $Z : [x,y,z]
\mapsto [\zeta x, \zeta^2 y, z]$ acts locally as a rotation
\[
  Z_{\phi}(t) = \beta t + O\big(t^2\big),
\]
where $\beta$ is characteristic of $Z$ and independent of $\phi$. Now,
on the one hand
\begin{gather*}
  Z_{\phi}\left(\phi([x,y,1])\right)^3
     = \phi\left(Z([x,y,1])\right)^3
    = \phi\big([\zeta x,\zeta^2 y,1]\big)^3
    = \alpha \zeta x + O\big(x^2\big),
\intertext{but also}
  Z_{\phi}\left(\phi([x,y,1])\right)^3
     = \left(\beta \phi([x,y,1]) + O\big(\phi^2\big)\right)^3
    = \beta^3 \phi([x,y,1])^3 + O\big(\phi^4\big)
    = \beta^3 \alpha x + O\big(x^2\big).
\end{gather*}
So $\beta^3 = \zeta$, or
\[
  \beta = \exp\left(\frac{2\pi\mi}{15} + \frac{2\pi\mi k}{3}\right)
\]
for some $k\in\{0,1,2\}$. But since $Z$ has order 5 we also know that
$\beta^5 = 1$, which in terms of $k$ means that
\[
  \frac{2\pi\mi}{3}+\frac{10\pi\mi k}{3}
  = \frac{2\pi\mi}{3}\left(1+5k\right) \in 2\pi\mi\mathbb{Z},
\]
or $\beta = \exp(\frac{4\pi\mi}{5})$ and at last we can conclude that
$Z$ corresponds to a rotation of $2\frac{2\pi\mi}{5}$ about $C$ in the
hyperbolic model.

\begin{figure}[t]
  \centering
  \includegraphics[scale=0.93]{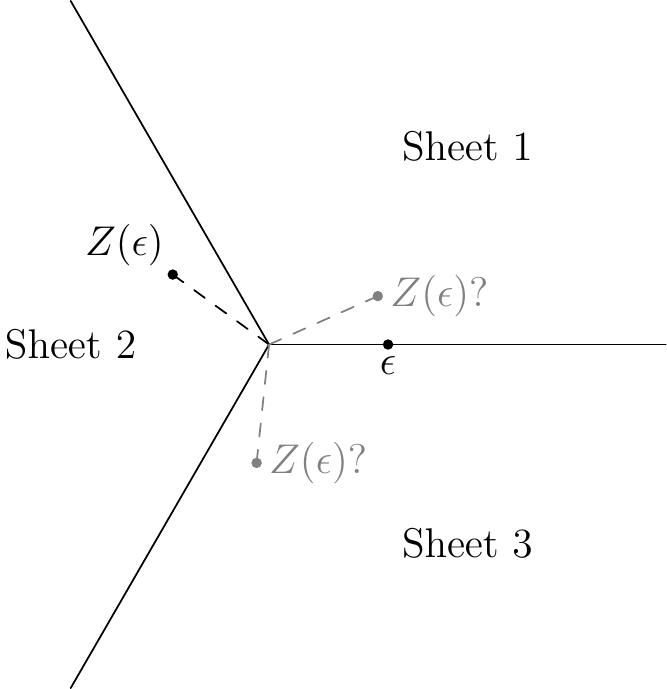}
  \caption{Intuitive action of $Z$ near $[0,0,1]$.}
  \label{fig:001Z}
\end{figure}

Intuitively we have unwrapped the three sheets coming together at
$[0,0,1]$ to obtain Fig.~\ref{fig:001Z} in $x$. We know that~$Z$
sends (say) $[\epsilon,y,1]$ to $[\zeta\epsilon,y',1]$ on some sheet
$y'$, which makes it one of the labelled destinations. But only one of
these gives an order 5 transformation so we know $Z$ completely.

\begin{figure}[t]
  \centering
  \includegraphics[scale=0.93]{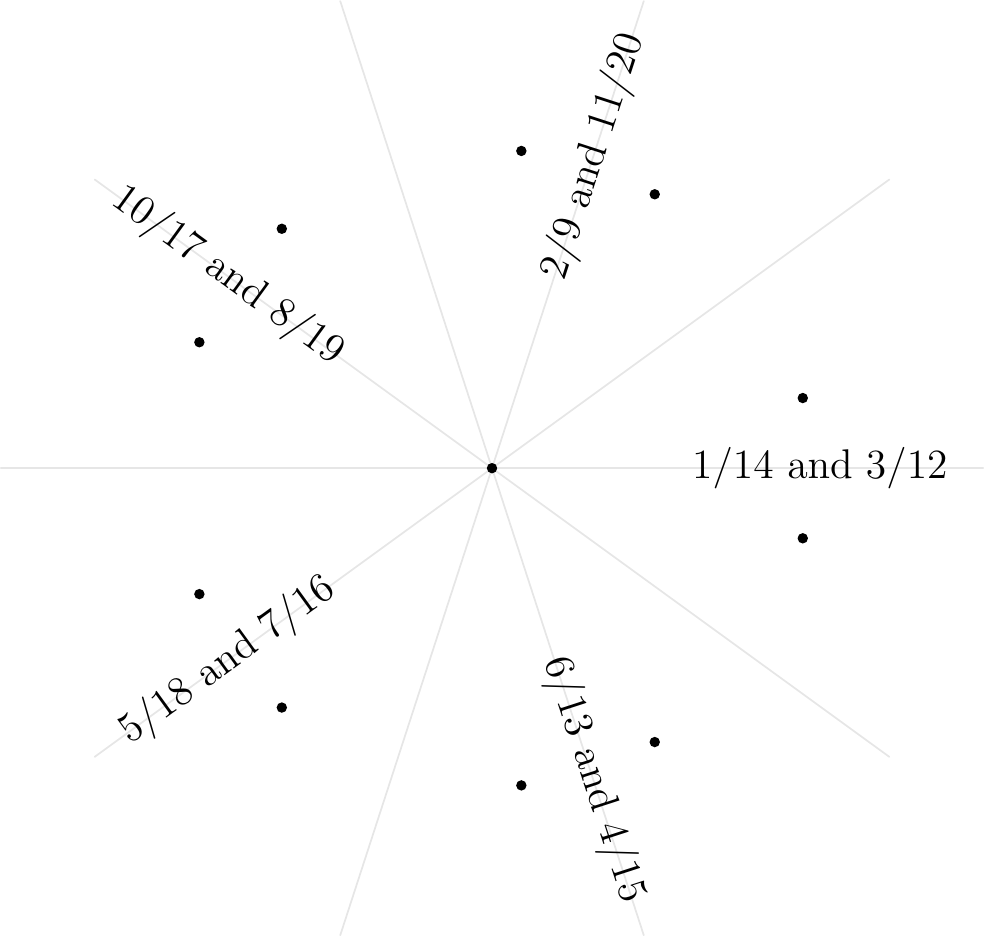}
  \caption{Hyperbolic polygonal edges in the Hulek--Craig model.}
  \label{fig:firstedges}
\end{figure}

Using this information, together with our knowledge that complex
conjugation on the HC model is the dashed ref\/lection in Fig.~\ref{fig:bringHyperbolic}, allows us to deduce the outline structure in
Fig.~\ref{fig:firstedges}. The dots are the branch-points of the HC
model and the grey lines are the images of the hyperbolic polygon's
edges under the isomorphism to the HC model. It remains to establish
which parts (and sheets) of each spoke in Fig.~\ref{fig:firstedges}
correspond to which hyperbolic edges (for example, does edge $1/14$
correspond to $x>0$ or $x < 0$, and what about $y$?).

Similar analysis of the other f\/ixed points of $Z$ will allow us
actually  to identify the remai\-ning~$P_i$. We f\/irst discover
\begin{itemize}\itemsep=0pt
\item Near $[0,1,0]$, $Z$ is a rotation of
  $3\left(\frac{2\pi}{5}\right)$.
\item Near $[1,0,0]_1 \sim [1,t^4,t]$, $Z$ is a rotation of
  $4\left(\frac{2\pi}{5}\right)$.
\item Near $[1,0,0]_2 \sim [1,t,t^4]$, $Z$ is a rotation of
  $\frac{2\pi}{5}$.
\end{itemize}
But hyperbolically, it is easy to see that a rotation of
$2\left(\frac{2\pi}{5}\right)$ about $C$ (which $Z$ is) is the same as
one of $4\left(\frac{2\pi}{5}\right)$ about $P_1$,
$3\left(\frac{2\pi}{5}\right)$ about $P_2$ or $\frac{2\pi}{5}$ about
$P_3$ so we can deduce that $[0,1,0] \leftrightarrow P_2$, $[1,0,0]_1
\sim [1,t^4,t] \leftrightarrow P_1$ and $[1,0,0]_2 \sim [1,t,t^4]
\leftrightarrow P_3$.

Therefore, edge $1/14$ corresponds to the real path from $[0,1,0]$ to
$[1,0,0]_1 \sim [1,t^4,t]$; referring to Fig.~\ref{fig:craigreal} we
see that this is the path where $y$ starts out large and positive near
$x=0$ (and remains positive). Edge~3 corresponds to the real path from
$[0,1,0]$ to $[1,0,0]_2\sim [1,t,t^4]$ which turns out to be the one
starting out large and negative near $x=0$ (and remaining negative).

The remaining paths ($y$ small near $x=0$) correspond to the diameter
of the hyperbolic model and have no large role to play in describing
the homology basis.

Other edges can now be obtained by applying a rotation of $2\pi/5$ on
the hyperbolic side and~$Z^3$ on the HC side. The results are in
Table~\ref{tab:hypcraigedges}.
\begin{table}[t]
\caption{Values for $[x,y,1]$ on hyperbolic edges.}
  \label{tab:hypcraigedges}
  \centering
  \begin{tabular}{rl|rl}
    1/14 & $[\mathbb{R}_+,\mathbb{R}_+,1]$ &
    2/9 & $[\zeta\mathbb{R}_+,\zeta^2\mathbb{R}_+,1]$ \\
    3/12 & $[\mathbb{R}_+,\mathbb{R}_-,1]$ &
    4/15 & $[\zeta^4\mathbb{R}_+,\zeta^3\mathbb{R}_-,1]$ \\
    5/18 & $[\zeta^3\mathbb{R}_+,\zeta\mathbb{R}_+,1]$ &
    6/13 & $[\zeta^4\mathbb{R}_+,\zeta^3\mathbb{R}_+,1]$ \\
    7/16 & $[\zeta^3\mathbb{R}_+,\zeta\mathbb{R}_-,1]$ &
    8/19 & $[\zeta^2\mathbb{R}_+,\zeta^4\mathbb{R}_-,1]$ \\
    10/17 & $[\zeta^2\mathbb{R}_+,\zeta^4\mathbb{R}_+,1]$ &
    11/20 & $[\zeta\mathbb{R}_+,\zeta^2\mathbb{R}_-,1]$
  \end{tabular}
\end{table}

\subsection{Riera and Rodr{\'{\i}}guez basis algebraically}\label{section4.4}

We are now in a position to express the Riera and Rodr{\'{\i}}guez basis on
this branched cover. Recall that{\samepage
\[
  \alpha_i = (2i-1) + (2i)
\]
as a prescription on which edges to traverse in the hyperbolic
model.}

This becomes a specif\/ication to look up the relevant edges in Table
\ref{tab:hypcraigedges}, and construct a path that has its main
component in the specif\/ied regions (circling $x=0$ and outside all
f\/inite branch points enough times to reach the correct sheets).  In
fact, just like Riera and Rodr{\'{\i}}guez we only need to construct
$\alpha_1$ and $\alpha_2$ and then repeatedly apply $(x,y) \mapsto
(\zeta x, \zeta^2 y)$ to obtain the rest.

\begin{figure}[t]
  \centering
  \includegraphics[scale=0.93]{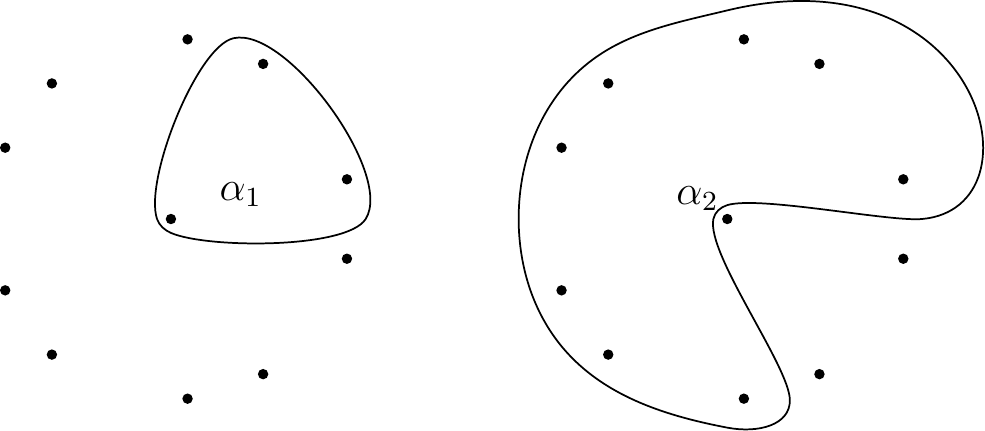}
  \caption{$\alpha_1$ and $\alpha_2$ homology cycles for the Hulek--Craig
  branched cover. Graphs of subsequent $\alpha_i$ are rotations of
  these by $\frac{2\pi\mi}{5}$.}
  \label{fig:alpha12}
\end{figure}

To be explicit and referring to Table \ref{tab:hypcraigedges},
$\alpha_1$ must go out along $x > 0$ with $y \gg 0$ near 0, loop
around inf\/inity until it can come back in to $x=0$ along a ray with
$\arg x = \frac{2\pi}{5}$ and $\arg y = \frac{4\pi}{5}$ before looping
around 0 until it can join up with the beginning again. A path
conforming to this description is shown in Fig.~\ref{fig:alpha12}.

Similarly $\alpha_2$ goes out along $x > 0$ with $y < 0$, loops and
comes back with argument of $x$ as~$-2\pi/5$ and argument of $y$ as
$6\pi/5$; it is also depicted in Fig.~\ref{fig:alpha12}.

Using the
software\footnote{Located at
\url{http://gitorious.org/riemanncycles}.} introduced in \cite{bn10} with
Klein's curve as an illustrative example, we may read these paths
into \texttt{extcurves} and convert them into a full basis
with the commands
\begin{verbatim}
> curve, hom, names := read_pic("homology.pic"):
> zeta := exp(2*Pi*I/5):
> trans := (x,y) -> [zeta^3*x,zeta*y]:
> for i from 1 to 3 do
    hom := [op(hom),
            transform_extpath(curve, hom[-2], trans),
            transform_extpath(curve, hom[-1], trans)];
  od:
\end{verbatim}

An immediate check to this calculation is provided by calculating the
intersection matrix of this constructed basis. The command
\begin{verbatim}
> Matrix(8, (i,j) -> isect(curve, hom[i], hom[j]));
\end{verbatim}
produces (with considerably less work and chance of error) precisely
the matrix claimed by Riera and Rodr{\'{\i}}guez, namely
\[
  \begin{pmatrix}
    0 & 1 & -1 & 1 & -1 & 0 & 1 & -1 \\
    -1 & 0 & 1 & -1 & 1 & 0 & 0 & 0 \\
    1 & -1 & 0 & 1 & -1 & 1 & -1 & 0 \\
    -1 & 1 & -1 & 0 & 1 & -1 & 1 & 0 \\
    1 & -1 & 1 & -1 & 0 & 1 & -1 & 1 \\
    0 & 0 & -1 & 1 & -1 & 0 & 1 & -1 \\
    -1 & 0 & 1 & -1 & 1 & -1 & 0 & 1 \\
    1 & 0 & 0 & 0 & -1 & 1 & -1 & 0
  \end{pmatrix}.
\]

Finally we can calculate the period matrix. This may be done analytically (as in \cite{riera}) or
numerically via the {\tt extcurves} package which calculates the
period matrix for any homology basis (implicitly using the Riemann period matrix given by {\tt algcurves[periodmatrix]}). Using the the
transformation from~\eqref{basistrans} both methods yield

\begin{theorem} The  homology cycles $\alpha_1$ and $\alpha_2$ for the Hulek--Craig
branched cover reproduce the cycles of Riera and Rodr{\'{\i}}guez and their corresponding
period matrix
\[
\tau = \tau_0
\begin{pmatrix}
  4  & 1  & -1 & 1  \\
  1  & 4  & 1  & -1 \\
  -1 & 1  & 4  & 1  \\
  1  & -1 & 1  & 4
\end{pmatrix},
\]
where $\tau_0$ is defined by \eqref{deftau0}.
\end{theorem}

We also note that the action on the homology basis $\alpha_{1,\ldots,8}$ associated
with the antiinvolution of the real structure is given by
\[
\mathcal{S}':=\begin{pmatrix} 0&0&0&0&0&0&-1&0\\
0 &0&0&0&0&-1&0&0\\  0&0&0&0&-1&0&0&0
\\  0&0&0&-1&0&0&0&0\\  0&0&-1&0&0&0
&0&0\\  0&-1&0&0&0&0&0&0\\  -1&0&0&0
&0&0&0&0\\  0&1&0&1&0&1&0&1\end{pmatrix}.
\]
It is algorithmic to show that
\[
\mathcal{S}:=T\mathcal{S}' T^{-1}=\begin{pmatrix} 1&0&0&0&1&0&0&0\\
0 &1&0&0&0&1&0&0\\  0&0&1&0&0&0&1&0
\\  0&0&0&1&0&0&0&1\\  0&0&0&0&-1&0&0
&0\\  0&0&0&0&0&-1&0&0\\  0&0&0&0&0&0
&-1&0\\  0&0&0&0&0&0&0&-1\end{pmatrix},
\]
where $T$ is the symplectic transformation (with respect to the canonical symplectic form)
\[
T= \begin{pmatrix} -1&2&1&2&-1&1&0&1
\\  0&1&0&-3&0&0&0&-1\\  -2&1&2&-1&0
&1&0&0\\  0&0&-1&0&-1&0&-1&0\\  2&-2
&-3&-2&0&-1&-1&-1\\  -1&0&1&3&0&1&0&1
\\  1&-1&-1&-2&0&-1&0&-1\\  -1&2&3&2
&0&1&1&1\end{pmatrix}.
\]
Here $\mathcal{S}$ is the canonical form for an antiholomorphic involution where there is one nondividing
real oval (see for example~\cite{vinnikov}), again showing there is one real oval.

\section{Vector of Riemann constants}\label{section5}

We shall now calculate the vector of Riemann constants for Bring's curve
determining various other quantities on the way. This vector together with
Riemann's theta function and the Abel map provide a
bridge between the analytic and algebraic structures of a Riemann
surface~$\mathcal{C}$, and as such are critical elements in the
implementation of the modern approach to integrable systems.

\subsection{The vector of Riemann constants}\label{section5.1}
Riemann established the fundamental result,
\[
\theta(e\,|\,\tau)=0 \ \Longleftrightarrow \
e\equiv
\mathcal{A}_Q\left(\sum_{i=1}\sp{g-1}P_i\right)-K_Q\in\Jac{\mathcal{C}},
\]
where $\theta$ is Riemann's theta function, $\mathcal{A}_Q$ is the
Abel map with base point~$Q\in\mathcal{C}$, $\tau$ is the period
matrix, $g$ is the genus of $\mathcal{C}$, $P_i\in\mathcal{C}$ and
the equivalence holds in the Jacobian
$\Jac\mathcal{C}=\mathbb{C}\sp{g}/(\mathbb{Z}\sp{g}+\tau\mathbb{Z}\sp{g})$.
(We are assuming that $g\ge1$ in what follows.) Both the period matrix $\tau$
and the Abel map depend on a choice of homology basis, the latter
through the basis of normalized holomorphic dif\/ferentials
$\boldsymbol{\omega}$ where
$\mathcal{A}_Q(P)=\int_Q\sp{P}\boldsymbol{\omega}$. The vector
$K_Q$ is known as the vector of Riemann constants\footnote{The
choice of sign of this vector depends on author. We will use that
of Fay \cite{fay} whose convention is the negative of Farkas and Kra \cite{farkaskra}.} (with
base point $Q$) and it also depends on the choice of homology
basis. Let ${\left\{\gamma_i\right\}_{i=1}\sp{2g}}=
\left\{\mathfrak{a}_i,\mathfrak{b}_i\right\}_{i=1}\sp{g}$ be our
choice of homology basis of $H_1(\mathcal{C},\mathbb{Z})$, where $
\mathfrak{a}_i$ and $\mathfrak{b}_i$ are canonically paired. One has that
\[
K_{Q j}  = -\frac{1}{2}\left( \tau_{jj} + 1\right) +
    \sum_{k\ne j}\oint_{\mathbf{\mathfrak{a}}_{k}}
    \omega_{k}(P)\int^{P}_{Q}\omega_{j}.
\]
Because of the integrations involved, both
$\tau$ and $K_Q$ are rather transcendental objects. One may also
express
\begin{equation}\label{vrcdelta}
K_Q\equiv\mathcal{A}_\ast\left(\Delta-(g-1)Q\right)
=\int_\ast\sp{\Delta}{\boldsymbol{\omega}}-(g-1)\int_\ast\sp{Q}{\boldsymbol{\omega}}
=\mathcal{A}_Q\left(\Delta\right),
\end{equation}
which holds for any base point $\ast$ of the Abel map and where
the degree $(g-1)$ divisor $\Delta$ is that of the Sz\"ego-kernel \cite{fay}. The
critical relation for us is the linear equivalence
\begin{equation*}\label{deltacanonical}
2\Delta\sim
    \mathcal{K}_\mathcal{C},
\end{equation*}
and so
\begin{equation}\label{Kcanonical}
2K_Q\equiv\mathcal{A}_Q(\mathcal{K}_\mathcal{C}).
\end{equation}
Here $\mathcal{K}_\mathcal{C}$ is the canonical divisor of
$\mathcal{C}$, the unique divisor class of degree $2g-2$ of any
meromorphic dif\/ferential on the curve, and (hereafter) $\sim$ denotes linear
equivalence.
Thus $\Delta$ gives a square root of the canonical
bundle or spin-structure on $\mathcal{C}$; the set $\Sigma$ of
divisor classes $\mathcal{D}$ such that
$2\mathcal{D}\sim\mathcal{K}_\mathcal{C}$ is called the set of
\emph{theta characteristics} of $\mathcal{C}$. The vector of Riemann constants
gives us the shift in the Jacobian necessary to identify spin structures
with the $2$-torsion points of the Jacobian. (Recall, an $N$-torsion point $x$ is such
that $Nx$ lies in the period lattice.)

\subsection{Symmetries and the vector of Riemann constants}\label{section5.2}

We now describe how symmetries may be used to restrict the vector
of Riemann constants by recalling some results we have established elsewhere.

Suppose that a curve has a nontrivial group of symmetries
$\Aut(\mathcal{C})$. Then the holomorphic dif\/ferentials of the curve,
$H\sp{1,0}(\mathcal{C},\mathbb{C})$, and
 the homology group $H_1(\mathcal{C},\mathbb{Z})$  are both
$\Aut(\mathcal{C})$-modules. Let $\sigma\in\Aut(\mathcal{C})$ and
denote the actions on these spaces by
\begin{equation*}
\sigma\sp\ast v_j= \sum_k v_k {L}\sp{k}_j,\qquad
\sigma_\ast\begin{pmatrix}\mathfrak{a}_i\\
  \mathfrak{b}_i\end{pmatrix}=M\begin{pmatrix}\mathfrak{a}_i\\
  \mathfrak{b}_i\end{pmatrix}:=
  \begin{pmatrix}{A}&B\\
    C&D\end{pmatrix}\begin{pmatrix}\mathfrak{a}_i\\
    \mathfrak{b}_i
  \end{pmatrix},
\end{equation*}
where ${L}\in GL(g,\mathbb{C})$, $M\in Sp(2g,\mathbb{Z})$ and $\{v_i\}$
is a (not necessarily normalized) basis of $H\sp{1,0}(\mathcal{C},\mathbb{C})$. Denote by $\Pi=\begin{pmatrix}\mathcal{A}\\
\mathcal{B}\end{pmatrix}$ the matrix of periods, where
$\mathcal{A}_{ij}=\int_{\mathfrak{a}_i}v_j$ and
$\mathcal{B}_{ij}=\int_{\mathfrak{b}_i}v_j$. Then
$\tau=\mathcal{B}\mathcal{A}\sp{-1}$ and $\boldsymbol{\omega}=
\boldsymbol{v}\mathcal{A}\sp{-1}$. The identity
$\oint_{\sigma_\ast\gamma}v=\oint_\gamma\sigma\sp\ast v$ (for any $\gamma\in H_1(\mathcal{C},\mathbb{Z})$)
 yields
the relation
\begin{equation}
\label{gensym}
M\Pi=\Pi {L}
\end{equation}
which restricts the period matrix $\tau$. Now using (\ref{vrcdelta})
and (\ref{Kcanonical}) we have that (with $\hat{L}=\mathcal{A} {L}\mathcal{A}\sp{-1}$, so as to
be working with normalized dif\/ferentials))
\[
2K_Q \hat{L}\equiv\int_\ast\sp{2\Delta}\sigma\sp\ast{\boldsymbol{\omega}}
-2(g-1)\int_\ast\sp{Q}\sigma\sp\ast{\boldsymbol{\omega}}
\] which yields
\[
2K_Q\big[\hat{L}-\id\big]\equiv\int_{2\Delta}\sp{\sigma(2\Delta)}{\boldsymbol{\omega}}-
2(g-1)\int_Q\sp{\sigma(Q)}{\boldsymbol{\omega}} .
\] If
$\mathcal{K}_\mathcal{C}$ is the divisor of a dif\/ferential $v$
then $\sigma\sp{-1}(\mathcal{K}_\mathcal{C})$ is the divisor of
$\sigma\sp\ast(v)$, whence the uniqueness of the canonical class
means that
$\sigma(\mathcal{K}_\mathcal{C})\sim\mathcal{K}_\mathcal{C}$ and
consequently $\sigma(2\Delta)\sim2\Delta$.
This shows that we have an action of $\Aut(\mathcal{C})$ on the theta characteristics and
\[
K_Q \hat{L}\equiv K_{\sigma(Q)}+\int_\Delta\sp{\sigma(\Delta)}{\boldsymbol{\omega}} ,
\]
the last integral also being a theta characteristic.
We have then the
identity on the Jacobian
\begin{equation*}
2K_Q\left[\hat{L}-\id\right]\equiv-
2(g-1)\int_Q\sp{\sigma(Q)}{\boldsymbol{\omega}}.
\end{equation*}
This then establishes
\begin{lemma}\label{lemtp} Suppose the automorphism $\sigma$ has order $N>1$. If $L-\id$ is invertible and $Q$ is a~fixed point of $\sigma$ then $K_Q$ is a $2N$-torsion point.
\end{lemma}

\begin{remark} We have the map
$\pi:\mathcal{C}\rightarrow\mathcal{C}/\langle \sigma\rangle$. Any holomorphic
dif\/ferential on $\mathcal{C}/\langle\sigma\rangle$ pulls back to an invariant
dif\/ferential on $\mathcal{C}$, and so the assumption that $L-\id$ is
invertible is equivalent to
$\mathcal{C}/\langle \sigma\rangle \cong\mathbb{P}\sp1$.
\end{remark}

\begin{corollary}\label{cortp} Assuming the conditions of  Lemma~{\rm \ref{lemtp}} and
that $\psi\in\Aut(\mathcal{C})$, then $\int_Q\sp{\psi(Q)}{\boldsymbol{\omega}}$ is a~$2N(g-1)$-torsion
point.
\end{corollary}

Although these simple results do not necessarily give
the best bound on the order of the torsion point for the vectors
involved, we see that, given a
suitable symmetry and f\/ixed point, we have that $K_Q$ is a torsion
point:
\begin{equation*}
2K_Q\mathcal{A}=\mathbf{n}\Pi\left[L-\id\right]\sp{-1}
=\mathbf{n}(M-\id)\sp{-1}\Pi=\mathbf{n}\left(\frac{1}{N}\sum_{k=1}\sp{N-1}
kM\sp{k}\right)\Pi.
\end{equation*}
The additional (we think) new idea we brought to this was to use some number theory
associated with $M$ to restrict the form of $K_Q$. Suppose
there exist $\mathbf{l}$, $\mathbf{m}\in\mathbb{Z}\sp{2g}$ such
that
\begin{equation}\label{sfdef}
\mathbf{m}=\mathbf{l}(M-\id),
\end{equation}
then
\[
\mathbf{m}\Pi=\mathbf{l}(M-\id)\Pi=
 \mathbf{l}\Pi\left[L-\id\right]
 \]
 and
\[
\left(2K_Q\mathcal{A}+\mathbf{l}\Pi\right)\left[{L}-\id\right]=(\mathbf{n}+\mathbf{m})\Pi  \in
\mathbb{C}\sp{g}.
\] The idea is to use the freedom in choosing
$\mathbf{l}$ here in~(\ref{sfdef}) to make $\mathbf{n}+\mathbf{m}$
as simple as possible; as we are only interested in $2K_Q$ modulo
the lattice we will have further restricted the choice of the
vector of Riemann constants. For example, if $\mathbf{m}$ could be
chosen arbitrarily then we could make
$\mathbf{n}+\mathbf{m}=\mathbf{0}$ and so $2K_Q$ would be a
lattice point. We implement the idea using the Smith normal
form of $M-\id$. Recall this means in the present context that we
may write
\[
M-\id=USV, \qquad S=\diag(d_1,\ldots,d_{2g}),\qquad d_i|d_{i+1},
\qquad U,V\in GL(2g,\mathbb{Z}).
\] The invertibility of $L-\id$ means
that $d_i\ge1$ and that (\ref{sfdef}) becomes
\[
\mathbf{m}V\sp{-1}=(\mathbf{l}U)S.
\]
Here we view $\mathbf{l}'=\mathbf{l}U$ as arbitrary and we are
interested in the constraints this places on $\mathbf{m}$. We have
$(\mathbf{m}V\sp{-1})_i=l_i'd_i$ and clearly the only constraints
arise for $d_i\ne1$. Given our earlier observation that here
$d_i\ge1$, we f\/ind that $\mathbf{m}$ is constrained only by
\[
\big(\mathbf{m}V\sp{-1}\big)_i\equiv0 \mod d_i,\qquad d_i>1.
\]
Thus, given a suitable symmetry and f\/ixed point~$Q$, the Smith
normal form of $M-\id$ enables us to restrict the possible torsion
points for~$2K_Q$. Considering further automorphisms and making
use of Corollary~\ref{cortp} may yield further restrictions. The f\/inal
step in evaluating $K_Q$ is the choice of the appropriate
half-period when taking the square root. This again may be
restricted by the symmetry but may also be decided numerically from
the $2^{2g}$ half-periods.

\subsection{Application to Bring's curve}\label{section5.3}

For Bring's curve we f\/ind that everything follows
from  study of the single (order 5) automorphism given in the
Hulek--Craig representation by
\[
  \phi : \  [\cx,\cy,\cz] \mapsto \big[\zeta^2\cx, \zeta^4 \cy, \cz\big].
\]
This has f\/ixed point $[0,0,1]$, or $Q=(0,0)$ in af\/f\/ine coordinates.

The f\/irst and easiest calculation is deriving its action on the
dif\/ferentials. We f\/ix the ordered basis of (unnormalized) holomorphic
dif\/ferentials
\begin{gather*}
  v_1  = \frac{(\cy^3-\cx)\dif \cx}
    {\partial_y\mathcal{C}(\cx,\cy,1)}, \qquad
  v_2  = \frac{(\cy^2\cx-1)\dif \cx}
    {\partial_y\mathcal{C}(\cx,\cy,1)}, \qquad
  v_3  = \frac{(\cy-\cx^2)\dif \cx}
    {\partial_y\mathcal{C}(\cx,\cy,1)}, \qquad
  v_4  = \frac{\cy(\cx^2-\cy)\dif \cx}
    {\partial_y\mathcal{C}(\cx,\cy,1)}.
\end{gather*}
The construction of such holomorphic
dif\/ferentials is algorithmic (see~\cite{DvH}).
It is easy to check that
\begin{gather*}
  \phi^*(v_1)  = \zeta v_1, \qquad
  \phi^*(v_2)  = \zeta^4 v_2,  \qquad
  \phi^*(v_3)  = \zeta^3 v_3, \qquad
  \phi^*(v_4)  = \zeta^2 v_4,
\end{gather*}
and so $\phi\sp\ast v_j= v_k {L}\sp{k}_j$
where
\begin{equation*}
 {L} =
  \begin{pmatrix}
    \zeta & 0 & 0 & 0 \\
    0 & \zeta^4 & 0 & 0 \\
    0 & 0 & \zeta^3 & 0 \\
    0 & 0 & 0 & \zeta^2
  \end{pmatrix}.
\end{equation*}
Thus there is no invariant dif\/ferential and ${L}-\id$ and so $\hat{L}-\id$
are invertible. With $Q=(0,0)$ we see the conditions of Lemma \ref{lemtp} are
satisf\/ied.

We note in passing that the dif\/ferential $v_3$ has a simple zero at $a=[0,0,1]$, a double
zero at $b=[0,1,0]$ and a triple triple zero at $c=[1,0,0]_2 \sim [1,t,t^4]$ for the required
total of $2g-2 = 6$. Thus we have $\mathcal{K}_\mathcal{C}\sim a+2b+3c$ expressing the
canonical divisor in terms of rational points of $\mathcal{C}$.

To proceed with our strategy of determining $K_Q$ we f\/irst
determine the action of $\phi$ on the homology cycles. With the program
\extcurves at
hand this is a simple computational matter, complicated only slightly by
the noncanonical nature of the paths we obtained in Section~\ref{section4.4}. We obtain $\phi_*(\gamma_i) =\sum_j M_{ij} \gamma_j$
where
\begin{equation*}
  M =
  \begin{pmatrix}
    0 & 0 & 0 & 1 & 0 & 0 & 0 & 0 \\
    -1 & 0 & 0 & -1 & 0 & 0 & 0 & 0 \\
    0 & -1 & 0 & 1 & 0 & 0 & 0 & 0 \\
    0 & 0 & -1 & -1 & 0 & 0 & 0 & 0 \\
    0 & 0 & 0 & 0 & -1 & 1 & -1 & 1 \\
    0 & 0 & 0 & 0 & -1 & 0 & 0 & 0 \\
    0 & 0 & 0 & 0 & 0 & -1 & 0 & 0 \\
    0 & 0 & 0 & 0 & 0 & 0 & -1 & 0
  \end{pmatrix}.
\end{equation*}
We remark that consideration of the equation (\ref{gensym}) for this
order f\/ive symmetry already imposes that
\[
  \mathcal{A}\sp{T}
  =
  \begin{pmatrix}
    a_1 & 0 & 0 & 0 \\
    0 & a_2 & 0 & 0 \\
    0 & 0 & a_3 & 0 \\
    0 & 0 & 0 & a_4
  \end{pmatrix}
  \begin{pmatrix}
    1 & -1-\zeta^4 & 1+\zeta^4+\zeta^3 & \zeta \\
    1 & -1-\zeta & 1+\zeta+\zeta^2 & \zeta^4 \\
    1 & -1-\zeta^2 & 1+\zeta^2+\zeta^4 & \zeta^3 \\
    1 & -1-\zeta^3 & 1+\zeta^3+\zeta & \zeta^2
  \end{pmatrix},
\]
for some unknown $a_i$. These $a_i$ are related by the remaining symmetries and
ultimately the period matrix~(\ref{RRstau}).

Continuing with our determination of $K_Q$, calculating the
Smith normal form of $M-1$ gives us unimodular
matrices $U$, $V$ such that
\begin{gather*}
  M - \id  = U\diag(1, 1, 1, 1, 1, 1, 5, 5) V
\end{gather*}
and our earlier constraint
$(\mathbf{m}V\sp{-1})_i\equiv0 \mod d_i$, $d_i>1$
takes the form
\begin{gather*}
  -m_5 +  m_6 - m_7 - 4m_8  \equiv 0 \pmod{5}, \\
  -m_1 + 2m_2 - 3m_3 - m_4 - 11m_5 + 6m_6 - m_7 - 34m_8  \equiv 0 \pmod{5}.
\end{gather*}
This gives 25 possible unique candidates for $2K_Q$; we can vary
$n_i$ arbitrarily (by adding an appropriate $\vec{l}$) without
essentially changing $2K_Q$ for (say) $i=2,3,4,6,7,8$ but then $n_1$
and $n_5$ are f\/ixed. Explicitly, every $2K_Q$ is equivalent to one
generated by
\[
  \vec{n} =
  \begin{pmatrix}
    n_1 & 0 & 0 & 0 & n_5 & 0 & 0 & 0
  \end{pmatrix}.
\]
By considering a further symmetry we f\/ind further that
$n_1 = n_5 = 3$ and at this stage we have
\begin{displaymath}
  2K_Q = \frac{1}{5}
  \begin{pmatrix}
    -12 & -3 & 3 & -3
  \end{pmatrix}
  + \tau_0
  \begin{pmatrix}
    -6 & -6 & 3 & 0
  \end{pmatrix}.
\end{displaymath}
To determining the appropriate 2-torsion point for square root of the canonical
bundle one could further study the action of the symmetries on the spin structures
or simply numerically test the vanishing of the theta function. The latter approach yields
that
\begin{theorem} For the Riera and Rodr{\'{\i}}guez homology basis of Bring's curve we have that the vector of
Riemann constants is
\begin{equation*}
    K_Q = \frac{1}{10}
    \begin{pmatrix}
      3 & 2 & -2 & -3
    \end{pmatrix}
    +
    \Im(\tau_0)
    \begin{pmatrix}
      1 & -2 & -2 & 1
    \end{pmatrix}\mi,
  \end{equation*}
where $\tau_0$ is defined by \eqref{deftau0}.
\end{theorem}

The transformation of theta characteristics
\[
\mathfrak{g}\cdot(\boldsymbol{a},\boldsymbol{b})=
(\boldsymbol{a},\boldsymbol{b})\mathfrak{g}^{-1}
+\frac12\big(\mathrm{diag}\big(CD^T\big),\mathrm{diag}\big(AB^T\big) \big)
\]
for any $\mathfrak{g}=\left(\begin{array}{cc}A&B\\C&D
\end{array}\right)\in\mathrm{Sp}(2g,\mathbb{Z})$ and characteristic
$(\boldsymbol{a},\boldsymbol{b})\in\mathbb{Q}^{2g}$ together with
the explicit representations of the symmetries yields
\begin{theorem}
Bring's curve has a unique invariant spin-structure.
\end{theorem}

\begin{remark} Klein's curve has a unique invariant spin structure~\cite{ks10}
and we have shown elsewhere that the vector of Riemann constants is the Abel image of this.
To show the analogous result for Bring's curve requires a better understanding of this
spin-structure.
\end{remark}

\subsection*{Acknowledgements}

We are grateful to Maurice Craig for helpful email exchanges and also to
an anonymous referee for careful reading and suggested improvements to the paper.

\pdfbookmark[1]{References}{ref}
\LastPageEnding


\begin{thebibliography}{99}
\footnotesize\itemsep=0pt

\bibitem{berndtiii}
Berndt B.C., Ramanujan's notebooks. {P}art~{III}, \href{http://dx.doi.org/10.1007/978-1-4612-0965-2}{Springer-Verlag}, New York,
  1991.

\bibitem{bn10}
Braden H.W., Northover T.P., Klein's curve, \href{http://dx.doi.org/10.1088/1751-8113/43/43/434009}{\textit{J.~Phys.~A: Math. Theor.}}
  \textbf{43} (2010), 434009, 17~pages, \href{http://arxiv.org/abs/0905.4202}{arXiv:0905.4202}.

\bibitem{breuer}
Breuer T., Characters and automorphism groups of compact {R}iemann surfaces,
  \textit{London Mathematical Society Lecture Note Series}, Vol.~280, Cambridge
  University Press, Cambridge, 2000.

\bibitem{begg}
Bujalance E., Etayo J.J., Gamboa J.M., Gromadzki G., Automorphism groups of
  compact bordered {K}lein surfaces. A~combinatorial approach, \textit{Lecture
  Notes in Mathematics}, Vol.~1439, Springer-Verlag, Berlin, 1990.

\bibitem{craig}
Craig M., A sextic {D}iophantine equation, \textit{Austral. Math. Soc. Gaz.}
  \textbf{29} (2002), 27--29.

\bibitem{craig04}
Craig M., On {K}lein's quartic curve, \textit{Austral. Math. Soc. Gaz.}
  \textbf{31} (2004), 115--120.

\bibitem{DvH}
Deconinck B., van Hoeij M., Computing {R}iemann matrices of algebraic curves,
  \href{http://dx.doi.org/10.1016/S0167-2789(01)00156-7}{\textit{Phys.~D}} \textbf{152/153} (2001), 28--46.

\bibitem{dye}
Dye R.H., A plane sextic curve of genus {$4$} with {$A_5$} for collineation
  group, \href{http://dx.doi.org/10.1112/jlms/52.1.97}{\textit{J.~London Math. Soc.}} \textbf{52} (1995), 97--110.

\bibitem{edge78}
Edge W.L., Bring's curve, \href{http://dx.doi.org/10.1112/jlms/s2-18.3.539}{\textit{J.~London Math. Soc.}} \textbf{18} (1978),
  539--545.

\bibitem{edge81}
Edge W.L., Tritangent planes of {B}ring's curve, \href{http://dx.doi.org/10.1112/jlms/s2-23.2.215}{\textit{J.~London Math. Soc.}}
  \textbf{23} (1981), 215--222.

\bibitem{farkaskra}
Farkas H.M., Kra I., Riemann surfaces, \textit{Graduate Texts in Mathematics},
  Vol.~71, Springer-Verlag, New York, 1980.

\bibitem{fay}
Fay J.D., Theta functions on {R}iemann surfaces, \textit{Lecture Notes in
  Mathematics}, Vol.~352, Springer-Verlag, Berlin, 1973.

\bibitem{hulek85}
Hulek K., Geometry of the {H}orrocks--{M}umford bundle, in Algebraic Geometry,
  {B}owdoin, 1985 ({B}runswick, {M}aine, 1985), \textit{Proc. Sympos. Pure
  Math.}, Vol.~46, Amer. Math. Soc., Providence, RI, 1987, 69--85.

\bibitem{ks10}
Kallel S., Sjerve D., Invariant spin structures on {R}iemann surfaces,
  \textit{Ann. Fac. Sci. Toulouse Math.~(6)} \textbf{19} (2010), 457--477,
  \href{http://arxiv.org/abs/math.GT/0610568}{math.GT/0610568}.

\bibitem{northover}
Northover T.P., Riemann surfaces with symmetry: algorithms and applications,
  Ph.D. thesis, Edinburgh University, 2011.

\bibitem{rauch}
Rauch H.E., Lewittes J., The {R}iemann surface of {K}lein with 168
  automorphisms, in Problems in Analysis (Papers Dedicated to {S}alomon
  {B}ochner, 1969), Princeton Univ. Press, Princeton, N.J., 1970, 297--308.

\bibitem{riera}
Riera G., Rodr{\'{\i}}guez R.E., The period matrix of {B}ring's curve,
  \textit{Pacific~J. Math.} \textbf{154} (1992), 179--200.

\bibitem{vinnikov}
Vinnikov V., Selfadjoint determinantal representations of real plane curves,
  \href{http://dx.doi.org/10.1007/BF01445115}{\textit{Math. Ann.}} \textbf{296} (1993), 453--479.

\bibitem{weber05}
Weber M., Kepler's small stellated dodecahedron as a {R}iemann surface,
  \href{http://dx.doi.org/10.2140/pjm.2005.220.167}{\textit{Pacific~J. Math.}} \textbf{220} (2005), 167--182.

\end{thebibliography}
\end{document}